\documentclass[12pt]{article}
\usepackage{theorem}
\usepackage{latexsym}
\usepackage{amssymb}
\usepackage{amsmath}
\usepackage{epic}
\def\Bx{\hfill{$\Box$}}

\newtheorem{thm}{Theorem}[section]

\newtheorem{cor}[thm]{Corollary}
\newtheorem{defn}[thm]{Definition}

\begin{document}

\renewcommand{\theequation}{\arabic{section}.\arabic{equation}}
\thispagestyle{empty}

\vskip 20pt
\begin{center}
{\bf REFINED RESTRICTED INVOLUTIONS}
\vskip 15pt
{\bf Emeric Deutsch}\\
{\it Department of Mathematics,}
{\it Polytechnic University,
Brooklyn, NY 11201}\\
{\tt deutsch@duke.poly.edu}

{\bf Aaron Robertson\footnote{
Homepage:  {\tt http://math.colgate.edu/$\sim$aaron/}}
}\\
{\it Department of Mathematics,}
{\it Colgate University,
Hamilton, NY 13346}\\
{\tt aaron@math.colgate.edu}

{\bf Dan Saracino}\\
{\it Department of Mathematics,}
{\it Colgate University,
Hamilton, NY 13346}\\
{\tt dsaracino@mail.colgate.edu}
\end{center}

\vskip 30pt
\begin{abstract}{\footnotesize \noindent
Define $I_n^k(\alpha)$ to be the set of involutions of
$\{1,2,\dots,n\}$ with exactly $k$ fixed points
which avoid the pattern $\alpha \in S_i$, for some $i \geq 2$, and
define $I_n^k(\emptyset;\alpha)$
to be the set of involutions of
$\{1,2,\dots,n\}$ with exactly $k$ fixed points
which contain the pattern $\alpha \in S_i$, for some $i \geq 2$,
exactly once.  Let
$i_n^k(\alpha)$ be the number of elements in $I_n^k(\alpha)$
and let $i_n^k(\emptyset;\alpha)$ be the number of elements in 
$I_n^k(\emptyset;\alpha)$.
We investigate $I_n^k(\alpha)$ and $I_n^k(\emptyset;\alpha)$ for all
$\alpha
\in S_3$.  In particular, we show that
$i_n^k(132)=i_n^k(213)=i_n^k(321)$,
$i_n^k(231)=i_n^k(312)$, $i_n^k(\emptyset;132)
=i_n^k(\emptyset;213)$, and $i_n^k(\emptyset;231)=i_n^k(\emptyset;312)$
for all $0 \leq k \leq n$. }
\end{abstract}
\vskip 50pt
\section*{\normalsize 1. Introduction}
\addtocounter{section}{+1}

Recall that $\pi \in S_n$
is called an involution if and only if $\pi^{-1}=\pi$.  Equivalently,
$\pi$ is an involution if and only if the cycle structure
of $\pi$ has no cycle of length longer than two.  In [RSZ],
the study of refined restricted permutations was initiated.
In order to describe the objects studied in [RSZ] and below
we have need of a few definitions.

Let $\pi \in S_n$ be a permutation
of
$\{1,2,\dots,n\}$ written in one-line notation.
Let $\alpha \in S_m$, $m \leq n$.
We say that $\pi$ {\it contains the pattern $\alpha$} if there exist
indices $i_1,i_2,\dots,i_m$ such that $\pi_{i_1} \pi_{i_2} \dots
\pi_{i_m}$ is equivalent to $\alpha$, where we define equivalence as
follows. Define $\overline{\pi}_{i_j}=
|\{x:\pi_{i_x} \leq \pi_{i_j}, 1 \leq x \leq m\}|$.  If
$\alpha = \overline{\pi}_{i_1} \overline{\pi}_{i_2} \dots
\overline{\pi}_{i_m}$ then
we say that $\alpha$ and $\pi_{i_1} \pi_{i_2} \dots \pi_{i_m}$ are
equivalent.  For example, if $\tau=124635$ then $\tau$ contains
the pattern $132$ by noting that $\tau_2 \tau_4 \tau_5 = 263$
is equivalent to $132$.  We say that $\pi$ is $\alpha$-{\it avoiding} if
$\pi$ does not contain the pattern $\alpha$. In our above example, $\tau$
is $321$-avoiding.

Let $S=\cup_{i \geq 2} S_i$.  Let $T$ be a subset of $S$
and $M$ be a multiset of $S$.
Define $S_n(T;M)$ to be 
the set of permutations in $S_n$ which avoid all patterns
in $T$ and contain each pattern in $M$ exactly once.  Let
$s_n(T;M)$ be the number of elements in $S_n(T;M)$. 
If $M = \emptyset$ we write $S_n(T)$ and $s_n(T)$.
Further, if $T$ or $M$ contain only one pattern, we omit the
set notation.

Consider the following refinement, introduced in [RSZ].
Define $S_n^k(T;M)$ to be the set of 
permutations in $S_n(T;M)$  with exactly
$k$ fixed points. Let
$s_n^k(T;M)$ be the number of elements in $S_n^k(T;M)$
where we omit $M$ and the set notation when appropriate.

In this paper, we are concerned with those permutations
in $S_n^k(T;M)$ which are involutions.  To this end, we define
$I_n^k(T;M)$ to be the set of involutions in $S_n^k(T;M)$
and we let $i_n^k(T;M)$ be the number of elements in
$I_n^k(T;M)$.  As before, we 
omit $M$ and the set notation when appropriate.

In [RSZ], it was shown that 
$s_n^k(132)=s_n^k(213)=s_n^k(321)$ and
$s_n^k(231)=s_n^k(312)$ for all $0 \leq k \leq n$.
In this paper we will show that the same result holds when restricting
our permutations to be involutions.  

The results $s_n^k(132)=s_n^k(321)$ and $i_n^k(132)= i_n^k(321)$
lend some evidence that there may be a restricted permutation
result concerning the cycle structure.  However,
for a given cycle structure $c$, in general,
the number of $132$-avoiding
permutations with cycle structure $c$ is {\it not } equal
 to the number of
$321$-avoiding permutations with cycle structure $c$.
As an example, consider $S_6(132)$ and $S_6(321)$.
(It should be noted that $n=6$ is the minimal $n$
such that the number of permutations
classified according to their cycle structure
differ by restriction.)
Below we give the permutations in each according
to their cycle structure.

$$
\begin{array}{lcccc}
\mathrm{Cycle \,\, structure}&|&S_6(132)&|&S_6(321)\\ \hline
1^6&|& 1 &|&1\\
1^42^1&|&5 &|&5\\
1^33^1&|&8 &|&8\\
1^22^2&|& 9&|&9\\
1^24^1&|&12 &|&12\\
1^12^13^1&|&20 &|&20\\
1^15^1&|&20 &|&20\\
2^3&|&5 &|&5\\
2^14^1&|&20 &|&18\\
3^2&|&8 &|&10\\
6^1&|&24 &|&24\\ \hline
\mathrm{Sum} &|&132&|&132
\end{array}
$$

Some results concerning restricted involutions along with
their fixed point refinement are known.  These are stated
in the following three theorems.  Other results
are given in [GM] and [GM2].

\begin{thm} (Simion and Schmidt, [SiS]\footnote{
Rodica Simion did not like
the SS acronym due to its unpleasant
connotation (see [Z]).  Hence, we use the
nonstandard SiS and hope that others will as well.})
Let $i_n(\alpha)$ be the number of $\alpha$-avoiding
involutions in $S_n$.   Let
$p_1 \in \{123,132,213,321\}$ and
$p_2 \in \{231,312\}$.
For $n \geq 1$,
$$
\begin{array}{l}
i_n(p_1)={n\choose {\lfloor\frac{n}{2}\rfloor}} \hskip 10pt \mathit{and}\\
 i_n(p_2)= 2^{n-1}.
\end{array}
$$

\end{thm}

\begin{thm} (Guibert and Mansour, [GM])
Let $i_n^k(132)$ be the number of $132$-avoiding
involutions in $S_n$ with $k$ fixed points.
For $0 \leq k \leq n$, 
$$
i_n^k(132) =
\left\{
\begin{array}{ll}
 \frac{k+1}{n+1} {n+1 \choose \frac{n-k}{2}}&\mathit{for \,\,}
k+n \,\, \mathit{even}\\ 
0&\mathit{for \,\,}
k+n \,\, \mathit{odd}.
\end{array}
\right.
$$ 
\end{thm}

\begin{thm} (Robertson, Saracino, and Zeilberger, [RSZ])
Let $\gamma \in S_n$ be given by $\gamma_i
=n+1-i$ for $1 \leq i \leq n$.  For $\pi \in S_n$, let $\pi^\star = \gamma
\pi
\gamma^{-1}$. Then, for all $\pi$, $\pi$ and $\pi^\star$ have the same
number of fixed points.  Furthermore, the number of occurrences
of the pattern $213$ (respectively $312$) in $\pi$ equals
the number of occurrences of the pattern $132$ (respectively $231$)
in $\pi^\star$.
\end{thm}

In the next section, we finish the enumeration of
$i_n^k(\alpha)$ for all $\alpha \in S_3$ and $0 \leq k \leq n$,
as well as provide some bijective results.
In the last section, we investigate
$I_n^k(\emptyset;\alpha)$ for all $\alpha \in S_3$ and
$0 \leq k \leq n$.

{\it Notation}  We note here that with a transposition
$(x y)$ we will always take $x<y$.

\section*{\normalsize 2.  Involutions Avoiding a Length Three
Pattern}
\addtocounter{section}{+1}
\setcounter{thm}{0}

In [SiS], Simion and Schmidt
completed the study of involutions avoiding a given pattern of length
three. Their results are given in Theorem 1.1 above.
As done in [RSZ], we refine the enumeration
problem by classifying restricted permutations according to the number of
fixed points.

We can see from the conjugation given in Theorem 1.3 that
$i_n^k(132)=i_n^k(213)$ and
$i_n^k(231)=i_n^k(312)$ for all $0 \leq k \leq n$.
In this section (Theorem  2.2) we show that
$i_n^k(321)=i_n^k(132)=i_n^k(213)$ for all $0 \leq k \leq n$
as well.

We note here that
since our permutations are involutions, we clearly
require $n+k$ to be even in all theorems below.

In the proofs below, we will use the following properties of
standard Young tableaux.  (For proofs of these
properties see [K], [K2], and [S].)  Let $\mathcal{Y}_\pi$
be the Young tableaux corresponding
(via the Robinson-Schensted algorithm) to $\pi \in S_n$.
Let $\mathcal{Y}_\pi$ have shape
$\lambda =
(\lambda_1 \geq \lambda_2 \geq \cdots )$,
 where $\lambda_i$ is the length of the $i^{\mathrm{th}}$ row.

\begin{enumerate}
\item $\lambda_1$ is the maximum length of an increasing subsequence
of $\pi$.
\item The length of the first column of $\mathcal{Y}_\pi$ is the
maximum length of a decreasing subsequence of $\pi$.
\item If $\pi$ is an involution, then the number of fixed points
of $\pi$ equals the number of odd length columns in $\mathcal{Y}_\pi$.
\item For $i \leq \frac{n}{2}$, the number of standard Young tableaux of
shape
$(n-i,i)$ (or its transpose via changing columns into rows)
is ${n \choose i} - {n \choose i-1}$.
\end{enumerate}

\begin{thm}  For $n \geq 1$,
$$
\begin{array}{rl}
i_n^0(123)=i_n^2(123)=&
\left\{
\begin{array}{ll}
{n-1 \choose \frac{n}{2}}&\mathit{for \,\, n \,\,
even}\\
0&\mathit{for \,\, n \,\,odd}\\
\end{array}
\right.
\\
\\
 i_n^1(123) =&
\left\{
\begin{array}{ll}
{n \choose \frac{n-1}{2}}&\mathit{for \,\, n \,\,odd}\\ 
0&\mathit{for \,\, n \,\,even}
\end{array}
\right.
\\
\\
i_n^k(123)=&0 \mathit{\,\,for \,\,} k \geq 3.
\end{array}
$$
\end{thm}

\noindent {\it Proof.}  Clearly for $k \geq 3$ we
have an occurrence of 123.  Hence, it remains to
prove the formulas for $k=0,1,2$.

Consider first $k=0$ so that $n$ is even.  Let
$\pi \in I_n^0(123)$.  Since $\pi$ is $123$-avoiding,
the longest increasing subsequence of $\pi$ has length
at most $2$.   Keeping in mind that $\pi$ has no
fixed
point, we use the above properties of standard
Young tableaux to see that
$$
i_n^0(123) = \sum_{j =0 \atop j \,\, even}^{\frac{n}{2}}
\left( {n \choose j} - {n \choose j-1} \right) = {n-1 \choose
\frac{n}{2}}.
$$

Next, consider $k=2$, so that again $n$ is even.
Since the total number of standard Young tableaux of two columns
on $\{1,2,\dots,n\}$ with $n$ even 
is ${n \choose \frac{n}{2}}$, we have
$$
i_n^2(123) = {n \choose \frac{n}{2}} - {n-1 \choose
\frac{n}{2}} = {n-1 \choose \frac{n}{2}-1}={n-1 \choose
\frac{n}{2}}.
$$

For $k=1$ we consider $i_{n-1}^1(123)$, which is equal to
the number of standard Young tableaux on $\{1,2,\dots,n\}$ with
at most $2$ columns, with $n$ odd (so that
exactly one of the columns is of odd length).
Hence,
$$
i_{n}^1(123) = \sum_{j=0}^{\frac{n-1}{2}}
\left( {n-1 \choose j} - {n \choose j-1} \right) = {n \choose
\frac{n-1}{2}}.
$$
\hfill{\Bx}

\noindent
{\it Remark.}  The case $i_n^1(123)$ also follows from Theorem
1.1 (originally done in [SiS]).

We now provide two bijections between $I_n^0(123)$ and
$I_n^2(123)$ since we see that they are enumerated
by the same sequence.

The first bijection uses standard Young tableaux.
For $\pi \in S_n$, denote by $SYT(\pi)$ the standard
Young tableau created by the Robinson-Schensted algorithm.
Let $SYT_n(2)$ be the set of all standard Young tableaux on
$n$ elements with at most $2$ columns with the lengths
of the columns having the same parity.
Now, let $\pi \in I_n^0(123)$ and consider $SYT(\pi)$.
From the properties of standard Young tableaux we see
that $SYT(\pi)$ has one or two columns, each of even
length.  Note that $n$ must be the bottom entry in
one of the columns.  Let $\gamma : SYT_n(2) \rightarrow SYT_n(2)$
be the map
which takes $n$ and places it on the bottom of the 
other column (even if empty).  For example,
$$
\gamma \left(
\begin{array}{cc}
1&3\\
2&4\\
5&\\
6
\end{array}
\right)
=
\left(
\begin{array}{cc}
1&3\\
2&4\\
5&6\\
\end{array}
\right).
$$

It is easy to check that for $\pi \in I_n^0(123)$,
 $\gamma(STY_n(\pi)) = STY_n(\tau)$ with $\tau \in I_n^2(123)$
and that $\gamma$ is a bijection.

The second bijection we present uses Dyck paths.  For
completeness we make the following definition.

\begin{defn}  Let $i \geq j \geq 0$ and $i+j \geq 2$ be even.  A partial
Dyck path is a path in
$\mathbb{R}^2$ from
$(0,0)$ to $(i,j)$ with $j > 0$ consisting of a sequence of steps
of length $\sqrt{2}$ and slope $\pm 1$ which does
not fall below the $x$-axis.   We denote
these two types of steps by $(1,1)$ and $(1,-1)$,
called up-steps and down-steps, respectively.  
If $j=0$ we call the path
a (standard) Dyck path. 
\end{defn}

{\bf Notation.} We will denote the set of
partial/standard Dyck paths from
$(0,0)$ to $(i,j)$ by $D(i,j)$.

We now describe, for completeness, a bijection from $S_n(123)$
to $D(2n,0)$ due to Krattenthaler [Kr]. 

Let $\mathcal{K}: S_n(123) \rightarrow D(2n,0)$ be the
bijection defined as follows.  Let $\pi_1 \pi_2 \cdots
\pi_n=\pi \in S_n(123)$.
Determine the
right-to-left maxima of $\pi$, i.e. $m=\pi_i$ is a 
right-to-left maximum if $m>\pi_j$ for all $j>i$.  
Let $\pi$ have right-to-left maxima $m_1<m_2<\dots<m_s$, so that
we may write
$$
\pi=w_sm_sw_{s-1}m_{s-1}\cdots w_1m_1,
$$
where the $w_i$'s are possibly empty.
Generate a Dyck path from $(0,0)$ to $(2n,0)$
as follows.
Read $\pi$ from right to left.  For each $m_i$ do $m_i-m_{i-1}$
up-steps (where we define $m_0=0)$.  For each $w_i$ do $|w_i|+1$
down-steps.

Using Krattenthaler's bijection, it is easy to check the following.

\begin{enumerate}
\item $\mathcal{K} |_{I_n(123)}$ produces a Dyck path that is symmetric
about the line
$x=n$.
\item $\mathcal{K} |_{I_n^0(123)}$ produces a Dyck path that has an even
number of peaks.
\item $\mathcal{K} |_{I_n^2(123)}$ produces a Dyck path that has an odd
number of peaks.
\end{enumerate}

For example, to prove 2 and 3, we note that for all $i$, $\pi_i$ is 
right-to-left maximum if and only if $i=1$, and that if there are two
fixed points then the righthand fixed point is a right-to-left maximum
but the lefthand fixed point is not.

Using facts 1--3, we define $\Gamma : I_n^0(123)
\rightarrow I_n^2(123)$ as follows.  Let 
$\pi \in I_n^0(123)$ and generate $\mathcal{K}(\pi)$,
which by the above properties must have a valley on
the line $x=n$, i.e. it must have a down-step
which ends on the line $x=n$
followed by an up-step.  To apply $\Gamma$,
turn the down-step into an up-step and the
up-step into a down-step.  
It is easy to check that $\Gamma$ is a bijection.

In the next theorem we find the surprising fact that
$i_n^k(132)=i_n^k(321)$ for all $0 \leq k \leq n$.

\begin{thm}  Let $\alpha \in \{132,213,321\}$.
For $0 \leq k \leq n$,
$$
i_n^k(\alpha) = \left\{
\begin{array}{ll}
 \frac{k+1}{n+1} {n+1 \choose \frac{n-k}{2}}&\mathit{for
\,\,} n+k \,\, \mathit{even}\\ 
 0&\mathit{for \,\,}
n+k \,\, \mathit{odd}.
\end{array}
\right.
$$ 
\end{thm}

\noindent {\it Proof.}  Due to Theorems 1.2 and 1.3, all that
remains is to prove the formula for the pattern $321$.
The proof for $321$ uses the properties of standard Young
tableaux.  Since $\pi \in I_n^k(321)$ may not contain
a decreasing subsequence of length greater than $2$,
we see that the Young tableaux corresponding to $\pi$ has shape
$(n-\frac{n-k}{2},\frac{n-k}{2})$.  Thus,
$$
i_n^k(321) = {n \choose \frac{n-k}{2}} - {n \choose \frac{n-k}{2}-1},
$$
which simplifies to the stated formula.
\hfill{\Bx}

We see, in particular, from Theorem 2.3, that the
$321$-avoiding derangement involutions of $\{1,2,\dots,2n\}$
and the $321$-avoiding involutions
of $\{1,2,\dots,2n-1\}$ with exactly one fixed point are
both enumerated by
$C_n=\frac{1}{n+1}{2n \choose n}$, the Catalan numbers.  To the best of
the authors' knowledge, these are new manifestation of the Catalan
numbers. Below, we provide a bijective
explanation of this fact, as a special case of the more general
bijection $\delta$ defined below.

It
is well-known that $|D(n,k)|= 
\frac{k+1}{n+1} {n+1 \choose \frac{n-k}{2}}$ (with $n+k$ even), the
formula given in Theorem 2.3.
Knowing this, we give a bijection from
$I_{n}^k(321)$ to $D(n,k)$. 
Note that $D(n,k)=\emptyset$ if $k<0$
or $k>n$.

Let $\pi \in I_{n}^k(321)$ with $n+k$ even and define the map $\delta:
I_{n}^k(321) \rightarrow D(n,k)$ as follows.
Write $\pi=\pi_1 \pi_2 \cdots \pi_{n}$.  If $\pi_i-i \geq 0$
then the $i^{\mathrm{th}}$ step in $\delta(\pi)$
is an up-step.   If $\pi_i-i<0$ then
the $i^{\mathrm{th}}$ step in $\delta(\pi)$
is a down-step.

We first show that $\delta(\pi) \in D(n,k)$ 
(i.e., that it does not fall below
the $x$-axis and that it ends at $(n,k)$).  Since $\pi$
is an involution, if we ignore all fixed points
in $\pi$, by the definition of
$\delta$, each down-step must be coupled with a remaining up-step to its
left.   Hence, for each $1 \leq i \leq n$,
$|\{j: \pi_j \geq j, j \leq i\}| \geq |\{j: \pi_j<j, j \leq i\}|$ thereby
showing that $\delta(\pi)$ does not fall
below the $x$-axis.  Since $k$ is the number of
fixed points in $\pi$ and we have $k$ more up-steps than
down-steps in $\delta(\pi)$, our ending height of $\delta(\pi)$ is clearly
$k$, thereby showing that $\delta(\pi) \in D(n,k)$.

To finish showing that $\delta$ is a bijection we 
provide $\delta^{-1}$.  Let $d \in D(n,k)$. 
Number the steps of $d$ from left
to right by $1,2,\dots,n$.
Proceeding from right to left across $d$, couple
each down-step with the closest uncoupled up-step to its 
left.  Take
the  two step numbers and create a transposition.
For the uncoupled up-steps (if any), take the step number
of and create a fixed point.
Once we have traversed $d$ we will have an
involution with $k$ fixed points. 

We now show that the resulting involution is $321$-avoiding.  We may
decompose $d$ as 
$$
u^{i_1} P u^{i_2} P u^{i_3} \cdots u^{i_{k}}P u^{i_{k+1}},
$$
with $k \geq 1$, $i_1,i_2,\dots,i_{k+1} \geq 0$, and
where $u^j$ stands for a sequence of $j$ consecutive up-steps
and the $P$'s are nonempty Dyck paths.  Hence, each occurrence of
$u$ in this decomposition is an uncoupled up-step and
yields a fixed point in $\delta^{-1}(d)$.  Furthermore, any
transposition in $\pi$ comes from an up-step and down-step
that both reside within the same $P$.

Note that a $321$ occurrence, if it exists,
may contain at most one fixed point.  Hence,
a $321$ occurrence must come from at least
two transposition, say $(ab)$ and $(cd)$.
Furthermore, from the description of $\delta^{-1}$,
we see that $(ab)$ and $(cd)$ must come
from the same $P$ in the decomposition given above
since if $(ab)$ comes from a $P$ to the left
of the $P$ from which $(cd)$ comes, then necessarily
$a,b<c,d$ and neither $c$ nor $d$ can be the smallest
element of the $321$ pattern.

First, let  
$z$ be a fixed point $(xy)$ and $(uv)$ be transpositions in
$\delta^{-1}(d)$, where $(xy)$ and $(uv)$ come from the same $P$ in the
decomposition above.
  From this decomposition
and the description of $\delta^{-1}$ we see that
either $z<x,y,u,v$ or $z>x,y,u,v$.  If $z<x,y,u,v$ then 
we have either a $123$ or a $132$ pattern.  If $z>x,y,u,v$ then
we have either a $123$ or a $213$ pattern.
Hence  a fixed point and at most two transpositions
cannot create a $321$ pattern.

We now let $(xy)$ and $(uv)$ be transpositions in $\delta^{-1}(d)$
which come from the same $P$ in the decomposition given above.
Without loss of generality, let $x<u$.  From the
description of $\delta^{-1}$ we must have $x<u<y<v$.
This ordering yields a $3412$ pattern, and thus no
occurrence of $321$.  Hence, any possible $321$ occurrence
must consist of one number from
each of three transpositions $(xy)$, $(uv)$, and $(wz)$.
We may assume that $x<u<w$ and conclude that
$x<u<w<y<v<z$.  This yields a $456123$ pattern, and thus
no occurrence of $321$.
 
An example is in order.
Consider $\pi=34125768 \in I_8^2(321)$.  Then $\delta(\pi)$
is the partial Dyck path shown below.

\vskip -30pt
\hskip 100pt
\setlength{\unitlength}{.75mm}
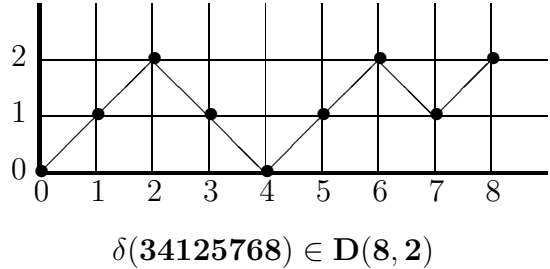
\begin{picture}(0,45)(-20,0)
\linethickness{.5mm}
\put(0,0){\line(0,1){30}}
\put(0,0){\line(1,0){90}}
\linethickness{.1mm}
\multiput(0,0)(0,10){3}{\line(1,0){90}}
\multiput(0,0)(10,0){9}{\line(0,1){30}}
\put(-5,-1){0}
\put(-5,9){1}
\put(-5,19){2}
\put(-1,-5){0}
\put(9,-5){1}
\put(19,-5){2}
\put(29,-5){3}
\put(39,-5){4}
\put(49,-5){5}
\put(59,-5){6}
\put(69,-5){7}
\put(79,-5){8}

\put(-1,-1){$\bullet$}
\put(9,9){$\bullet$}
\put(19,19){$\bullet$}
\put(29,9){$\bullet$}
\put(39,-1){$\bullet$}
\put(49,9){$\bullet$}
\put(59,19){$\bullet$}
\put(69,9){$\bullet$}
\put(79,19){$\bullet$}

\linethickness{.5mm}
\put(0,0){\line(1,1){20}}
\put(20,20){\line(1,-1){20}}
\put(40,0){\line(1,1){20}}
\put(60,20){\line(1,-1){10}}
\put(70,10){\line(1,1){10}}
\end{picture}
\vskip 10pt
\begin{center}
{$\mathbf{\delta(34125768) \in D(8,2)}$}
\end{center}

For the inverse,
we traverse the above partial Dyck path from
right to left to get the involution (in cycle
notation) $(8) (6 \, 7) (5) (2 \, 4) (1 \, 3) =
34125768$.

We may also use a bijection to $D(n,k)$ to offer
an alternative proof for the patterns $132$ and
$213$ (which by Theorem 1.3 are essentially the same).
Let $\pi \in I_{n}^k(213)$ with $n+k$ even and consider the bijection
$\zeta: I_n^k(213) \rightarrow D(n,k)$ defined as follows. 

Create two columns, the left column designated the
{\it up column}, denoted $UC$, and the right column
designated the {\it down column}, denoted $DC$.
Let $\pi=\pi_1 \pi_2 \cdots \pi_{n}$.  Read $\pi$ from
left to right while performing the following algorithm.

I.  If $\pi_i=i$, move to the next row, place $i$
in $UC$, and move down another row.

II.  If $\pi_i>i$, let $x$ be the largest entry in $DC$'s
row.  If $x$ does not exist, set $x=0$.

\hskip 10pt a) If $\pi_i>x$, place $\pi_i$ in $DC$ and $i$
in $UC$.
\vskip -10pt
\hskip 10pt b) If $\pi_i<x$, move to the next row and
place $\pi_i$ in $DC$ and $i$ in $UC$.  Furthermore,
move
\vskip -10pt
\hskip 10pt
 any element $y \in DC$, $y<\pi_i$, that is in a
row above $\pi_i$ to the row in which $\pi_i$ was placed.

After placing all elements of $\pi$ into $UC$ or
$DC$, compute $(u_1,u_2,\dots,u_t)$ and
$(d_1,d_2,\dots,d_t)$, where $u_i$ is the number of entries
in the $i^{\mathrm{th}}$ row of $UC$ and $d_i$
is the number of entries
in the $i^{\mathrm{th}}$ row of $DC$.
(Note that some of the $d_i$'s may be $0$.)  
The partial Dyck path given by
$u^{u_1} d^{d_1} u^{u_2} d^{d_2} \cdots u^{u_t} d^{d_t}$,
where $u^j$ is $j$ consecutive up-steps and
$d^j$ is $j$ consecutive down-steps, is $\zeta(\pi)$

To show that
$\zeta$ is a bijection, we give $\zeta^{-1}$.
Let $d \in D(n,k)$.  Traversing $d$ from left to
right label the up-steps in order (starting with $1$).
Once this is done, traversing $d$ from right to left,
label the down-steps in order starting with
the next number (one more than the number of up-steps
in $d$).    

Call an up-step and a down-step to the right of
the up-step {\it matching} if
the line segment connecting their midpoints does
not intersect the partial Dyck path.
Using the labeling of steps given above,
create a transposition of the labels for every pair of matching
up-steps and down-steps.  If an up-step has no
matching down-step, create a fixed point with its label.

We now provide a sketch that the resulting permutation is $213$-avoiding.
We may
decompose $d$ as 
$$
u^{i_1} P u^{i_2} P u^{i_3} \cdots u^{i_{k}}P u^{i_{k+1}},
$$
with $k \geq 1$, $i_1,i_2,\dots,i_{k+1} \geq 0$, and
where $u^j$ stands for a sequence of $j$ consecutive up-steps
and the $P$'s are nonempty Dyck paths.  Hence, each occurrence of
$u$ in this decomposition is an unmatched up-step and
yields a fixed point in $\delta^{-1}(d)$.  Furthermore, any
transposition in $\delta^{-1}(d)$ comes from an up-step and down-step
that both reside within the same $P$.

Let $f$ be a fixed point $\delta^{-1}(d)$.  Note that
all elements to the left of $f$ in $\delta^{-1}(d)$
are either fixed points or are larger than $f$.
It follows that a $213$ occurrence cannot contain
two fixed points.  Also, if $x$ and $y$ are not fixed
points, then only $xfy$ may be a $213$ pattern.
However, this implies that $(xy)$ is not
a transposition, i.e., that $(xa)$ and $(yb)$ are
the transpositions (with $a \neq y$).  Further,
$(xa)$ must come from a $P$ in the decomposition to the
left of the up-step corresponding to $f$ and
$(yb)$ must come from a $P$ in the decomposition to the
right of the up-step corresponding to $f$.  But this
implies that $x>y,b$ and so $xfy$ is not an
occurrence of $213$.
Thus, a $213$ occurrence cannot contain a fixed point.

Now assume that $(xy)$ and $(uv)$, $x<u$, create a $213$ pattern.
The only ordering which yields a $213$ pattern is $x<y<u<v$. 
However, this is not possible since we have an up-step
(corresponding to $u$) with a higher label than a
down-step (corresponding to $y$).

The last remaining case to consider is 
$(xy), (uv)$, $(wz)$, $x<u<w$, creating a $213$ pattern.
We must have the ordering $x<u<w<v<y<z$ in order
to have a $213$ pattern (in fact, two such patterns).  However, such an
ordering is not possible since
the path matching $u$ and $v$ will intersect
one of the paths matching $x$ and $y$ or $w$ and $z$.

To illustrate $\zeta$, consider the following example.
Let $\pi=689751423 \in I_{9}^1(213)$.  We find that our
up and down columns are 
$$
\begin{array}{ccc}
\mathrm{UC}&\vline&\mathrm{DC}\\
\hline
1,2,3&\vline&8,9\\
4&\vline&6,7\\
5&\vline&\\
\end{array}.
$$

From here we get $(u_1,u_2,u_3)=(3,1,1)$ and
$(d_1,d_2,d_3)=(2,2,0)$.
Hence, $\zeta(\pi)$
is the partial Dyck path given below (ignoring the labels
and dotted lines).

\hskip 80pt
\setlength{\unitlength}{.75mm}
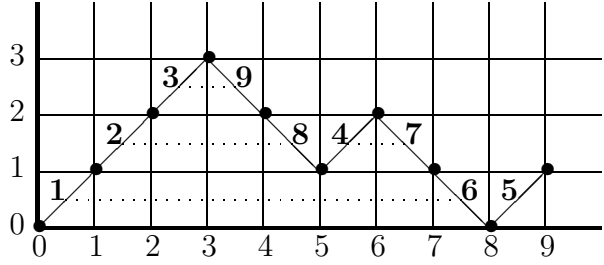
\begin{picture}(0,45)(-20,0)
\linethickness{.5mm}
\put(0,0){\line(0,1){40}}
\put(0,0){\line(1,0){100}}
\linethickness{.1mm}
\multiput(0,0)(0,10){4}{\line(1,0){100}}
\multiput(0,0)(10,0){10}{\line(0,1){40}}
\put(-5,-1){0}
\put(-5,9){1}
\put(-5,19){2}
\put(-5,29){3}
\put(-1,-5){0}
\put(9,-5){1}
\put(19,-5){2}
\put(29,-5){3}
\put(39,-5){4}
\put(49,-5){5}
\put(59,-5){6}
\put(69,-5){7}
\put(79,-5){8}
\put(89,-5){9}

\put(-1,-1){$\bullet$}
\put(9,9){$\bullet$}
\put(19,19){$\bullet$}
\put(29,29){$\bullet$}
\put(39,19){$\bullet$}
\put(49,9){$\bullet$}
\put(59,19){$\bullet$}
\put(69,9){$\bullet$}
\put(79,-1){$\bullet$}
\put(89,9){$\bullet$}

\linethickness{.5mm}
\put(0,0){\line(1,1){30}}
\put(30,30){\line(1,-1){20}}
\put(50,10){\line(1,1){10}}
\put(60,20){\line(1,-1){20}}
\put(80,0){\line(1,1){10}}

\put(2,5){{\bf 1}}
\put(12,15){{\bf 2}}
\put(22,25){{\bf 3}}
\put(35,25){{\bf 9}}
\put(45,15){{\bf 8}}
\put(52,15){{\bf 4}}
\put(65,15){{\bf 7}}
\put(75,5){{\bf 6}}
\put(82,5){{\bf 5}}

\linethickness{.1mm}
\dottedline{2}(5,5)(75,5)
\dottedline{2}(15,15)(45,15)
\dottedline{2}(25,25)(35,25)
\dottedline{2}(55,15)(65,15)

\end{picture}
\vskip 10pt
\begin{center}
{$\mathbf{\zeta(689751423) \in D(9,1)}$}
\end{center}
\vskip 10pt

For the inverse, 
note that a dotted line connects 
an up-step with its matching down-step (if it exists).
Using this information and the labels 
on the above partial Dyck path we can immediately construct
(in cycle
notation) $(1 \,6) (2 \,8)(3 \,9) (4 \, 7) (5) =
689751423$.

The next theorem finishes this section.

\begin{thm} Let $\alpha \in \{231,312\}$.  For $n \geq 1$
and $0 \leq k \leq n$,
$$
i_n^k(\alpha)=
\left\{
\begin{array}{ll}
2^{\frac{n-k-2}{2}}\left( {\frac{n+k}{2} \choose 
\frac{n-k}{2}} +
{{\frac{n+k-2}{2}} \choose 
\frac{n-k}{2}} \right)&\mathit{for \,\,}n+k \,\, \mathit{even}\\ 
0&\mathit{for \,\,}n+k \,\, \mathit{odd.}\\
\end{array}
\right.
$$
\end{thm}

\noindent {\it Proof.} In [SiS] it is remarked that
$S_n(\{231,312\}) = I_n(231)$.  Hence,
 $S_n^k(\{231,312\}) = I_n^k(231)$ for all
$0 \leq k \leq n$. 
This last equality,
coupled with Theorem 2.8 in [MR], gives the stated formula.
\hfill{\Bx}

\section*{\normalsize 3.  Involutions Containing a Length Three
Pattern Exactly Once}
\addtocounter{section}{+1}
\setcounter{thm}{0}

We can see from the conjugation given
in Theorem 1.3 that
$i_n^k(\emptyset;132)=i_n^k(\emptyset;213)$ and
$i_n^k(\emptyset;231)=i_n^k(\emptyset;312)$ for all $0 \leq k \leq n$.
In this section, we show that these are the only
equalities for patterns of length three.
We note again that
since our permutations are involutions, we clearly
require $n+k$ to be even in all theorems below.

\begin{thm}  For $n \geq 1$,
$$
\begin{array}{ll}
i_n^3(\emptyset;123) = \frac{3}{n}{n \choose \frac{n-3}{2}}&
\mathit{for\,\,}n \geq 3 \,\, \mathit{odd, \,\, and}\\ \\
i_n^k(\emptyset;123) =0&\mathit{otherwise}.
\end{array}
$$
\end{thm}

\noindent {\it Proof.}  We start with some
at-first-sight unrelated results.

Recall that $D(i,j)$ is the set of partial/standard
Dyck paths from $(0,0)$ to $(i,j)$.  Define $d(n,j)$ to be
the size of $D(2n-j-1,j-1)$ for $ j \geq 0$.  Since a
step ending at $(2n-j-1,j-1)$ is either a down-step
from $(2n-j-2,j)$ or an up-step from $(2n-j-2,j-2)$ we see
that
\begin{equation}
d(n,j)=d(n,j+1)+d(n-1,j-1).
\end{equation}

By definition we have
$d(n,1)=C_{n-1}$.  From (3.1) we  get
$d(n,2)=C_{n-1}$ as well.  Rearranging
(3.1) and making the change of variables $j \mapsto j+1$
and $n \mapsto n+1$ we
get
\begin{equation}
d(n,j)=d(n+1,j+1)-d(n+1,j+2).
\end{equation}

From (3.2) we have
$$
\sum_{j=2}^n d(n,j) = \sum_{j=2}^n
\left(d(n+1,j+1)-d(n+1,j+2)\right)=d(n+1,3).
$$

Applying (3.1) again we see that
\begin{equation}
\sum_{j=2}^n d(n,j) = d(n+1,3)=d(n+1,2)-d(n,1) = C_{n}-C_{n-1}.
\end{equation}

Now consider $\mathcal{K}: S_n(123) \rightarrow D(2n,0)$, 
Krattenthaler's  bijection as described in section 2.
 
Let $\pi$ have right-to-left maxima $m_1<m_2<\dots<m_s$, so that
we may write
$$
\pi=w_sm_sw_{s-1}m_{s-1}\cdots w_1m_1,
$$
where the $w_i$'s are possibly empty.

We notice that if $\pi_j=n$ then since $m_s=n$
we have $|w_s|=j-1$.  We consider the
adumbrated permutation (which is technically not a permutation,
but obviously corresponds uniquely to a
permutation of the same length)
$$
\pi^\star=m_sw_{s-1}m_{s-1}\cdots w_1m_1,
$$
(i.e. $\pi$ with $w_s$ removed).  Using
the algorithmic steps of $\mathcal{K}$,
we may abuse notation and write  $\mathcal{K}(\pi^\star)$
to mean $\mathcal{K}(\pi)$
with its last step removed.  This partial Dyck path is in  
$D(2n-j-1,j-1)$.  To see this, note that $\mathcal{K}(\pi^\star)$
ends at $(2n-j,j)$ but the last step {\it must}
be an up-step.  Hence, the
number of permutations in
$S_n(123)$ with $\pi(j)=n$ is $d(n,j)$ for any $1 \leq j \leq n$.

At last we turn our attention to $I_n^k(\emptyset;123)$.
We first argue that if $k \neq 3$ then $i_n^k(\emptyset;123)=0$.
Clearly, if $k>3$ we have more than one occurrence of $123$.  Hence,
we assume $k<3$.
Let $\pi \in I_n^k(\emptyset;123)$ with $k <3$ and let our $123$
pattern be the subsequence $a b c$ in $\pi$.
It is easy to see that if we let
$\pi=\pi(1)b\pi(2)$ then $\pi(1)$ is a permutation
of $\{a,b+1,b+2,\dots,c-1,c+1,\dots,n\}$ and $\pi(2)$ is
a permutation of $\{1,2,\dots,a-1,a+1,\dots, b-1,c\}$.
Since $\pi$ is an involution, we see that we must
have both $a$ and $c$ as fixed points.
This in turn implies that $b$ must be fixed,
since $b$ is preceded by $b-1$ entries,
contradicting our assumption that $k<3$.

Thus, we restrict our attention to $k=3$, whereby
our $3$ fixed points create the single
$123$ occurrence.  Call these fixed points $a<b<c$.
From above we see that we
must have $b = \frac{n+1}{2}$ and $n$ odd in order
for $b$ to be a fixed point.

Since we are restricted to involutions, the placement of 
$1,2,\dots, b-1,c$ completely defines $\pi \in I_n^3(\emptyset;123)$.
Thus, $1,2,\dots, b-1,c$ must be $123$-avoiding
and can be identified  uniquely with some $\tau \in S_b(123)$
with $\tau(j)=b$, where $j=j'-(b-1)$ and $j'$ is defined by $\pi^{-1}(c)$.
Since $\tau(1)=a$ we must have $\tau(1)\neq c$ so that
$j \neq 1$.  Since the number of permutations
in $S_n(123)$ with $\pi(j)=n$ is $d(n,j)$
and $b=\frac{n+1}{2}$ we have, using (3.3),
\begin{equation}
i_n^3(\emptyset;123) = \sum_{j=2}^{\frac{n+1}{2}} d\left(\frac{n+1}{2},j
\right)
= C_{\frac{n+1}{2}} - C_{\frac{n-1}{2}},
\end{equation}
which simplifies to the stated formula.
\hfill{\Bx}

As a consequence of Theorem 3.1, we obtain the
following obvious corollary.

\begin{cor} For $n \geq 3$, $i_n(\emptyset;123)= \frac{3}{n}{n \choose
\frac{n-3}{2}}$.
\end{cor}

The next theorem for the pattern $132$ was first proved in
[GM].  This combined with Theorem 1.3 yields the following
theorem, which we include for completeness.

\begin{thm}  For $n \geq 3$, $0 \leq k \leq n$, and $\alpha \in
\{132,213\}$,
$$
i_n^k(\emptyset;\alpha) =
\left\{
\begin{array}{ll}
 \frac{k+1}{n-1}{n-1 \choose \frac{n+k}{2}}&
\mathit{for \,\,}n+k \,\, \mathit{even \,\, and \,\,} k \neq 0\\ \\
0&\mathit{otherwise}.
\end{array}
\right.
$$

\end{thm}

Summing $i_n^k(\emptyset;\alpha)$ for $\alpha \in \{132,213\}$ over
$k$ gives us the following result, first given in [GM]. 

\begin{cor} For $n \geq 3$ and
$\alpha \in \{132,213\}$, $i_n(\emptyset;\alpha)= {n-2 \choose
\lfloor\frac{n-3}{2}\rfloor}$.
\end{cor}

\begin{thm}  For $n \geq 4$, $0 \leq k \leq n$, and $\alpha \in
\{231,312\}$,
$$
i_n^k(\emptyset;\alpha) =
\left\{
\begin{array}{ll}
(k-1)2^{\frac{n-k-6}{2}}\left( {\frac{n+k}{2}-2 \choose 
\frac{n-k}{2}-1} +
2{{\frac{n+k}{2}-3} \choose 
\frac{n-k}{2}-1} 
+{\frac{n+k}{2}-4 \choose 
\frac{n-k}{2}-1} \right)&\mathit{for \,\,}n+k \,\, \mathit{even}\\ \\
0&\mathit{for \,\,}n+k \,\, \mathit{odd.}\\
\end{array}
\right.
$$
\end{thm}

\noindent {\it Proof.}  Let $a(n,k)=i_n^k(\emptyset;231)$
and $b(n,k) = i_n^k(231)$ for $0 \leq k \leq n$.
Let $\pi \in I_n^k(\emptyset;231)$.  Write
$\pi=\pi(1)n\pi(2)j$; if $j=1$ then
$\pi(1)=\emptyset$ and if $j=n$ then $\pi(2) =\emptyset$. 

For $j=n$ we clearly have $\pi(1) \in I_{n-1}^{k-1}(\emptyset;231)$.
For $j<n$ we consider two cases:  the $231$ pattern is to the left of $n$
and the $231$ pattern is to the right of $n$.  We argue that the
$231$ pattern cannot include $n$.
To see this, assume otherwise and let $ynx$
be the $231$ pattern.  If $x \neq j$ we must have $j>y$
so that $x<j$ and $(xy)$ is not a transposition of $\pi$.
This implies that $\pi^{-1}(x)nx$ and $ynx$
are distinct $231$ patterns (since $(xy)$ is not
a transposition of $\pi$), a contradiction.  If $x=j$ then we have $j<y$.
Since $(jn)$ is a transposition of $\pi$, we know
that $\pi^{-1}(y) \neq j$.  Hence, 
$ynj$ and $yn\pi^{-1}(y)$ are two distinct occurrences
of $231$, again a contradiction.

First, consider the case where $\pi(1)$ contains the pattern.
Note that we must have $\pi(2) = (n-1)(n-2) \cdots (j+2)(j+1)$
and that $\pi(1) \in I_{j-1}^{k-1} \cup I_{j-1}^k$,
depending upon the parity of $n+j$.  Hence, this case
contributes
$$
\sum_{j =1 \atop n+j \,\, even}^{n-2} a(j-1,k-1) + \sum_{j =1 \atop n+j
\,\, odd}^{n-1} a(j-1,k).
$$

Next, consider the case where $\pi(2)j$ contains the pattern.
In this case it is easy to see that $\pi(2)=(n-2)  (n-1)$
and that $j=n-3$.  Thus,
this case contributes $b(n-4,k-2)$ to the total.

Summing over all $j$ we get
$$
a(n,k)=a(n-1,k-1) + b(n-4,k-2) + 
\sum_{{j =1} \atop {n+j \,\, even}}^{n-2} a(j-1,k-1) + \sum_{{j =1} \atop
{n+j\,\, odd}}^{n-1} a(j-1,k)
$$

which, using $b(n,k)=2b(n-2,k)+b(n-1,k-1)$ given in
[MR], yields
\begin{equation}
a(n,k)=2a(n-2,k)+a(n-1,k-1)+b(n-6,k-2)+b(n-5,k-3).
\end{equation}

Define the generating functions
 $A_k(x) = \sum_{n \geq 0} a(n,k)x^n$ and
$B_k(x) =
\sum_{n
\geq 0} b(n,k)x^n$.  Using $b(n,k)=2b(n-2,k)+b(n-1,k-1)$ it is
easy to show that
\begin{equation}
B_k(x) = \frac{x^k(1-x^2)}{(1-2x^2)^{k+1}}.
\end{equation}
Since $b(n,k)=s_n^k(231,312)$ (shown in the
proof of Theorem 2.4) we have
(from Theorem 2.9 in [MR])
\begin{equation}
B_k(x)= \sum_{{n \geq 1} \atop {n+k \,\, even}}
2^{\frac{n-k-2}{2}}\left( {\frac{n+k}{2}
\choose 
\frac{n-k}{2}} +
{{\frac{n+k-2}{2}} \choose 
\frac{n-k}{2}} \right) x^n.
\end{equation}

From (3.5) we have $A_k(x) = 2x^2A_k(x) +x A_{k-1}(x)
+x^6B_{k-2}(x) + x^5B_{k-3}(x)$.  Using (3.5) and (3.6) we get
$$
A_k(x) = \frac{(k-1)x^{k+2}(1-x^2)^2}{(1-2x^2)^k}
= (k-1)x^3(1-x^2)B_{k-1}(x).
$$

To obtain the stated formula for $a(n,k)$, we extract the coefficient of
$x^n$ in
$A_k(x)$ using the above equation and (3.7) and simplify.
\hfill{\Bx}

Summing $i_n^k(\emptyset;\alpha)$ for $\alpha \in \{231,312\}$ over
$k$ gives us the following nice formula.

\begin{cor} For $n \geq 5$ and
$\alpha \in \{231,312\}$, $i_n(\emptyset;\alpha)= (n-1)2^{n-6}$.
\end{cor}

\noindent
{\it Remark.}  For $n=4,6,8,\dots$, $i_n(\emptyset;\alpha)
= i_{2n-4}^2 (\emptyset;\alpha)$ for $\alpha \in \{231,312\}$.

The last remaining pattern to consider in this section is $321$.

\begin{thm}  For $n \geq 3$, $0 \leq k \leq n$, 
$$
i_n^k(\emptyset;321) = \frac{k(k+3)}{n+1} {n+1 \choose \frac{n-k}{2}-1}.
$$

\end{thm}

\noindent{\it Proof.}  
Let $i(n,k) = i_n^k(\emptyset;321)$.
We first show that
$$
i(n,k) = \sum_{f=1 \atop f \,\, odd}^{n-k}
i(n-f,k-1) C_{\frac{f-1}{2}}
 + \sum_{f=2 \atop f \,\, even}^{n-k}  
i_{n-f}^k(321) C_{\frac{f}{2}},
$$
which is equivalent, by Theorem 2.3, to 
\begin{equation}
i(n,k) = \sum_{f=1 \atop f \,\, odd}^{n-k}i(n-f,k-1)  C_{\frac{f-1}{2}}
+ \sum_{f=2 \atop f \,\, even}^{n-k}  \frac{k+1}{n-f+1}
{n-f+1 \choose \frac{n-k-f}{2}} C_{\frac{f}{2}},
\end{equation}
where we have initial conditions  $i(n,0)=0$ for all $n \geq 3$.

To see that $i(n,0)=0$ for all $n$ 
let $cba$ be the $321$ pattern in $\pi$. In order to avoid another $321$
pattern, to the left (right) of $b$ we cannot have an element larger
(smaller) than $b$,  except $c$ ($a$) itself. Hence, $b$ is a fixed point.
Thus, the restriction of having  exactly one $321$ pattern implies a fixed
point must be present (see Theorem 6.4 in [RSZ] for further details).
 Hence, we may let $f$ be the smallest fixed
point in $\pi \in I_n^k(\emptyset;321)$.  We separate
the argument into two cases:  $f$ odd and $f$ even.

First, let $f$ be odd and write $\pi = \pi(1)f\pi(2)$.  In order
for $\pi$ to contain exactly one occurrence of $321$ we
must have $\pi(1) \in I_{f-1}^0(321)$ and
$\pi(2) \in I_{n-f}^{k-1}(\emptyset;321)$.
To see that we require $\pi(1)\in I_{f-1}^0(321)$
assume otherwise, that is that $\pi(1)$ is not
an involution.  Since $\pi$ is an
involution and $f$ is odd, there exist $x \neq y$ both in $\pi(1)$ with
$x,y >f$.  This produces two occurrences of $321$:
$xf \pi(x)$ and $yf\pi(y)$, a contradiction.
(As an aside, this shows that an odd fixed point
cannot be part of a $321$ occurrence.)
Next, since $\pi(1)\in I_{f-1}^0(321)$, we necessarily
must have $\pi(2) \in I_{n-f}^{k-1}(\emptyset;321)$. 
Summing over valid $f$ and using Theorem 2.3
(for $k=0$) we get
$$
\sum_{f=1 \atop f \,\, odd}^{n-k}
i(n-f,k-1) C_{\frac{f-1}{2}}
$$
in this case.

Next, consider $f$ even.  Again, write $\pi = \pi(1)f\pi(2)$.
Since $\pi$ is an involution and $f$ is even we must
have $x \in \pi(1)$ with $x>f$.  This gives the $321$
occurrence $xf \pi(x)$.  Thus, only one such $x$ may exist.
Furthermore, $\pi(2)$ must be $321$-avoiding.  Now, consider
the $f$ leftmost entries in $\pi$:  $\tau = \tau(1)x\tau(2)f$.
Note that $\tau$ is a $321$-avoiding permutation on $f$ elements.
Furthermore, $\tau$ does not contain the element $\pi^{-1}(x)$.
Thus, $\tau$ is a permutation of $\{1,2,\dots,\pi^{-1}(x)-1,\pi^{-1}(x)+1,
\dots,f-1,f,x\}$.  By letting $i \in \tau$ become $i-1$
if $\pi^{-1}(x)+1 \leq i \leq f$ and letting $x$ become $f$ we obtain
$\tau^\star \in I_f^0(321)$.  Next consider $\pi(2)^\star =x\pi(2)$, i.e.
$\pi(2)$ with $x$ in the first position. 
As before, $\pi(2)^\star$ may be identified with $\sigma \in
I_{n-f+1}^{k-1}(321)$ with the added condition that
$\sigma(1) \neq 1$ since we know that $\pi(2)^\star (1) = x > 
\pi^{-1}(x)$.  Next, we  have
that the number of $321$-avoiding involutions
of $\{1,2,\dots,n-f+1\}$ with $k-1$
fixed points and $1$ not a fixed point is
$i_{n-f+1}^{k-1}(321) - i_{n-f}^{k-2}(321)$
(where $ i_{n-f}^{k-2}(321)$ counts the number of
such permutations with $1$ being a fixed point).
Noting that $i_{n}^{k-1}(321)-i_{n-1}^{k-2}(321)
= i_{n-1}^k(321)$ and
summing over valid $f$ 
we get
$$
\sum_{f=2 \atop f \,\, even}^{n-k}  
i_{n-f}^{k}(321) C_{\frac{f}{2}}.
$$
Combining the two cases' results proves (3.8).

We must now show that (3.8) along with the initial
conditions yields $i(n,k) = \frac{k(k+3)}{n+1} {n+1 \choose
\frac{n-k}{2}-1}$.

We first show that
$$
\frac{k+3}{n+1} {n+1 \choose
\frac{n-k}{2}-1}
=
\sum_{f=2 \atop f \,\, even}^{n-k}  \frac{k+1}{n-f+1}
{n-f+1 \choose \frac{n-k-f}{2}} C_{\frac{f}{2}},
$$
i.e., that 
\begin{equation}
\frac{k+3}{n+1} {n+1 \choose
\frac{n-k}{2}-1}
=
\sum_{i=1}^{\frac{n-k}{2}}  \frac{k+1}{n-2i+1}
{n-2i+1 \choose \frac{n-k-2i}{2}} C_{i}.
\end{equation}

For $0 \leq k \leq n$, denote the lefthand side of (3.9) by $f(n,k)$ and
the righthand side of (3.9) by $g(n,k)$.  
  It is straightforward
to show that for $k \geq 1$,
$f(n,k)=f(n-1,k+1)+f(n-1,k-1)$ and
$g(n,k)=g(n-1,k+1)+g(n-1,k-1)$, where we
define $f(n,k)=0$ and $g(n,k)=0$ if $n<k$.   Since
$f(2,2)=g(2,2)$, to prove that (3.9) holds
it is sufficient to show that $f(n,0)=g(n,0)$ for
all $n \geq 2$.

By Theorem 3.1, we see that $f(n,0)= \frac{3}{n+1}
{n+1 \choose \frac{n}{2}-1}=i_n^3(\emptyset;123)$.
From (3.4), this gives us $f(n,0) = C_{\frac{n}{2}+1} - C_{\frac{n}{2}}$,
where $C_n$ is the Catalan number.
Next, since 
$$
\begin{array}{ll}
g(n,0) &= \sum_{i=1}^{\frac{n}{2}}  \frac{1}{n-2i+1}
{n-2i+1 \choose \frac{n}{2}-i} C_{i} \\ \\
& = \sum_{i=1}^{\frac{n}{2}} C_{\frac{n}{2}-i} C_i\\ \\
& = \sum_{i=0}^{\frac{n}{2}} C_{\frac{n}{2}-i} C_i - C_{\frac{n}{2}}\\ \\
& =  C_{\frac{n}{2}+1} - C_{\frac{n}{2}}\\ \\
& = f(n,0)
\end{array}
$$
we have proven (3.9).

We now have
\begin{equation}
i(n,k) = \sum_{f=1 \atop f \,\, odd}^{n-k}i(n-f,k-1)  C_{\frac{f-1}{2}}
+ \frac{k+3}{n+1} {n+1 \choose
\frac{n-k}{2}-1},
\end{equation}
with initial conditions  $i(n,0)=0$ for all $n \geq 2$.

We use this and induction on $n+k$ to prove that
$i(n,k) = \frac{k(k+3)}{n+1} {n+1 \choose
\frac{n-k}{2}-1}$.  Since this holds
for $i(1,1)$ and $i(2,0)$, we may assume
that $i(n-f,k-1) = \frac{(k-1)(k+2)}{n-f+1} {n-f+1 \choose
\frac{n-f-k-1}{2}}$.  Substitution into (3.10) gives
$$
i(n,k) = \sum_{i=1}^{\frac{n-k}{2}}
\frac{(k-1)(k+2)}{n-2i+2} {n-2i+2 \choose
\frac{n-k}{2}-i}  C_{i-1}
+ \frac{k+3}{n+1} {n+1 \choose
\frac{n-k}{2}-1}.
$$
Hence, we must show that
\begin{equation}
\sum_{i=1}^{\frac{n-k}{2}}
\frac{k+2}{n-2i+2} {n-2i+2 \choose
\frac{n-k}{2}-i}  C_{i-1}
=
\frac{k+3}{n+1} {n+1 \choose
\frac{n-k}{2}-1}.
\end{equation}
Denote by $h(n,k)$ the lefthand side of (3.11) and keep
$f(n,k)$ as the notation for the righthand side of (3.11).
It is straightforward to show that $h(n,k)=h(n-1,k+1)+h(n-1,k-1)$
and that $h(1,1)=f(1,1)$ and $h(2,2)=f(2,2)$.
To prove (3.11), it is sufficient to show that $h(n,0)=f(n,0)$
for all $n \geq 2$.  Since
$$
\begin{array}{ll}
h(n,0) &= \sum_{i=1}^{\frac{n}{2}}
\frac{2}{n-2i+2} {n-2i+2 \choose
\frac{n}{2}-i}  C_{i-1}\\ \\
&= \sum_{i=0}^{\frac{n}{2}-1}
\frac{2}{n-2i} {n-2i \choose
\frac{n}{2}-i-1}  C_{i}\\ \\
&=\sum_{i=0}^{\frac{n}{2}-1} C_{\frac{n}{2}-i} C_i \\ \\
& = \sum_{i=0}^{\frac{n}{2}} C_{\frac{n}{2}-i} C_i - C_{\frac{n}{2}} \\ \\
& =  C_{\frac{n}{2}+1}  - C_{\frac{n}{2}} \\ \\
& = f(n,0)
\end{array}
$$
we have proven (3.11), thereby proving the theorem.
\hfill{$\Box$}

From the proof of Theorem 3.7 we obtain
Corollary 3.9 below, for which we have need of the
following definition.

\begin{defn}  Let $dp(n,k) \in D(n,k)$ and let $dp_x(n)$
be a Dyck path with $2n$ steps starting at $(x,0)$.
For $1 \leq i \leq \frac{n-k}{2}$,
we call a lattice path which results
from
$dp(n-2i,k) \cup dp_{n-2i}(i)$
a modified Dyck path with a single drop from height $k$,
and denote the set of all such modified
Dyck paths by $MDP(n;k)$.
\end{defn}

Using this definition, we can give the following, the
proof of which is a direct consequence of (3.9). 

\begin{cor} For $n \geq 2$ and $0 \leq k \leq n$
with $n+k$ even,
$|MDP(n;k)|=\frac{k+3}{n+1} {n+1 \choose
\frac{n-k}{2}-1}$.
\end{cor}

Comparing Corollary 3.9 with the number of
partial Dyck paths, we find that
$|MDP(n;k)|=|D(n,k+2)|$.
We explain this via a bijection.

Let $pdp(n-2i,k) \circ dp(i)$ be the decomposition of
an element in $MDP(n;k)$ where $pdp$ stands
for partial Dyck path and $dp$ stands for
(standard) Dyck path.  To obtain an
element in $D(n,k+2)$ we perform the following
steps.

Concatenate one up-step to the end of $pdp(n-2i,k)$.
To the end of this new up-step concatenate $d(i)$
and remove the last step of $d(i)$ (necessarily a
down-step).  The result is an element of
$D(n,k+2)$.

For the inverse, perform the following steps
to $pdp(n,k+2) \in D(n,k+2)$.  Add a down-step
to the end of $pdp(n,k+2)$.  Next, traverse $pdp(n,k+2)$
from left to right and locate the last occurrence of
two consecutive up-steps whose second step has
ending point on the line $y=k+2$.
From these two up-steps, remove the up-step closest to the origin.  We now
have a partial Dyck path ending at height $k$ and a Dyck
path lying $k+1$ units above the $x$-axis.  Move the
Dyck path left 1 unit and down $k+1$ units.  The result is
a member of $MDP(n;k)$.

We illustrate this bijection with an example.
Consider the following member of $MDP(10;2)$.
\vskip -20pt
\hskip 70pt
\setlength{\unitlength}{.75mm}
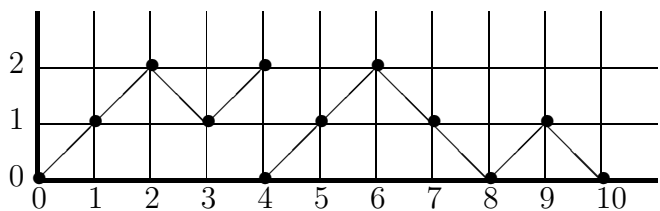
\begin{picture}(0,45)(-20,0)
\linethickness{.5mm}
\put(0,0){\line(0,1){30}}
\put(0,0){\line(1,0){110}}
\linethickness{.1mm}
\multiput(0,0)(0,10){3}{\line(1,0){110}}
\multiput(0,0)(10,0){11}{\line(0,1){30}}
\put(-5,-1){0}
\put(-5,9){1}
\put(-5,19){2}
\put(-1,-5){0}
\put(9,-5){1}
\put(19,-5){2}
\put(29,-5){3}
\put(39,-5){4}
\put(49,-5){5}
\put(59,-5){6}
\put(69,-5){7}
\put(79,-5){8}
\put(89,-5){9}
\put(99,-5){10}

\put(-1,-1){$\bullet$}
\put(9,9){$\bullet$}
\put(19,19){$\bullet$}
\put(29,9){$\bullet$}
\put(39,-1){$\bullet$}
\put(39,19){$\bullet$}
\put(49,9){$\bullet$}
\put(59,19){$\bullet$}
\put(69,9){$\bullet$}
\put(79,-1){$\bullet$}
\put(89,9){$\bullet$}
\put(99,-1){$\bullet$}

\linethickness{.5mm}
\put(0,0){\line(1,1){10}}
\put(10,10){\line(1,1){10}}
\put(20,20){\line(1,-1){10}}
\put(30,10){\line(1,1){10}}
\put(40,0){\line(1,1){10}}
\put(50,10){\line(1,1){10}}
\put(60,20){\line(1,-1){10}}
\put(70,10){\line(1,-1){10}}
\put(80,0){\line(1,1){10}}
\put(90,10){\line(1,-1){10}}
\end{picture}
\vskip 20pt
We add an up-step to the end of the partial
Dyck path  and remove the
last step of the modified Dyck path to get
the following lattice path.

\hskip 70pt
\setlength{\unitlength}{.75mm}
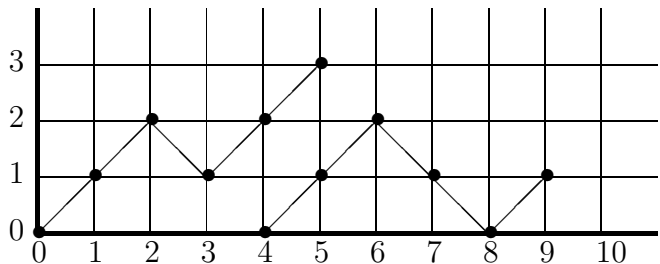
\begin{picture}(0,45)(-20,0)
\linethickness{.5mm}
\put(0,0){\line(0,1){40}}
\put(0,0){\line(1,0){110}}
\linethickness{.1mm}
\multiput(0,0)(0,10){4}{\line(1,0){110}}
\multiput(0,0)(10,0){11}{\line(0,1){40}}
\put(-5,-1){0}
\put(-5,9){1}
\put(-5,19){2}
\put(-5,29){3}
\put(-1,-5){0}
\put(9,-5){1}
\put(19,-5){2}
\put(29,-5){3}
\put(39,-5){4}
\put(49,-5){5}
\put(59,-5){6}
\put(69,-5){7}
\put(79,-5){8}
\put(89,-5){9}
\put(99,-5){10}

\put(-1,-1){$\bullet$}
\put(9,9){$\bullet$}
\put(19,19){$\bullet$}
\put(29,9){$\bullet$}
\put(39,-1){$\bullet$}
\put(39,19){$\bullet$}
\put(49,9){$\bullet$}
\put(49,29){$\bullet$}
\put(59,19){$\bullet$}
\put(69,9){$\bullet$}
\put(79,-1){$\bullet$}
\put(89,9){$\bullet$}

\linethickness{.5mm}
\put(0,0){\line(1,1){10}}
\put(10,10){\line(1,1){10}}
\put(20,20){\line(1,-1){10}}
\put(30,10){\line(1,1){10}}
\put(40,0){\line(1,1){10}}
\put(40,20){\line(1,1){10}}
\put(50,10){\line(1,1){10}}
\put(60,20){\line(1,-1){10}}
\put(70,10){\line(1,-1){10}}
\put(80,0){\line(1,1){10}}

\end{picture}
\vskip 20pt
To create an element of $D(n,k+2)$ we concatenate the
Dyck path with its last step removed to the end of
the partial Dyck path and get the following.

\vspace*{60pt}
\hskip 70pt
\setlength{\unitlength}{.75mm}
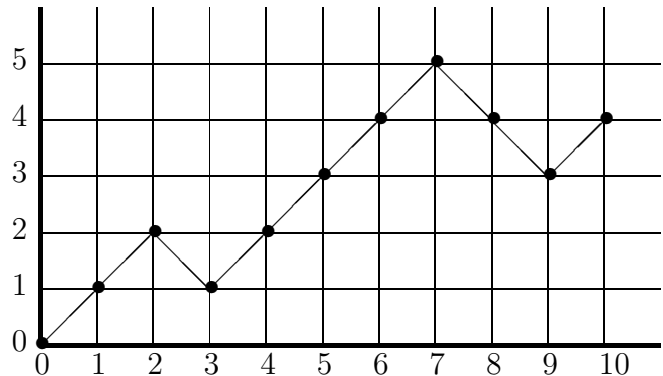
\begin{picture}(0,45)(-20,0)
\linethickness{.5mm}
\put(0,0){\line(0,1){60}}
\put(0,0){\line(1,0){110}}
\linethickness{.1mm}
\multiput(0,0)(0,10){6}{\line(1,0){110}}
\multiput(0,0)(10,0){11}{\line(0,1){60}}
\put(-5,-1){0}
\put(-5,9){1}
\put(-5,19){2}
\put(-5,29){3}
\put(-5,39){4}
\put(-5,49){5}
\put(-1,-5){0}
\put(9,-5){1}
\put(19,-5){2}
\put(29,-5){3}
\put(39,-5){4}
\put(49,-5){5}
\put(59,-5){6}
\put(69,-5){7}
\put(79,-5){8}
\put(89,-5){9}
\put(99,-5){10}

\put(-1,-1){$\bullet$}
\put(9,9){$\bullet$}
\put(19,19){$\bullet$}
\put(29,9){$\bullet$}
\put(49,29){$\bullet$}
\put(39,19){$\bullet$}
\put(59,39){$\bullet$}
\put(69,49){$\bullet$}
\put(79,39){$\bullet$}
\put(89,29){$\bullet$}
\put(99,39){$\bullet$}

\linethickness{.5mm}
\put(0,0){\line(1,1){10}}
\put(10,10){\line(1,1){10}}
\put(20,20){\line(1,-1){10}}
\put(30,10){\line(1,1){10}}
\put(50,30){\line(1,1){10}}
\put(40,20){\line(1,1){10}}
\put(60,40){\line(1,1){10}}
\put(70,50){\line(1,-1){10}}
\put(80,40){\line(1,-1){10}}
\put(90,30){\line(1,1){10}}

\end{picture}

\section*{\normalsize References}
\footnotesize

[GM] O. Guibert and T. Mansour,
Restricted $132$-involutions and Chebyshev
polynomials,\\ arXiv:math.CO/0201136 (2002).

[GM2] O. Guibert and T. Mansour,
Some statistics on restricted $132$ involutions,
arXiv:math.CO/0206169 (2002).

[K] D. E. Knuth, \underline{The Art of Computer Programming}, vol. 3,
Addison-Wesley, Reading, MA, 1973.

[K2] D. E. Knuth, Permutations, matrices, and generalized
Young tableaux, {\it Pacif. J. Math.} {\bf 34} (1970), 709-727. 

[Kr] C. Krattenthaler, Permutations with restricted
patterns and Dyck paths, {\it Advances in Applied
Math.} {\bf 27} (2001), 510-530.

[MR] T. Mansour and A. Robertson, Refined restricted permutations
avoiding subsets of patterns of length three,
arXiV:math.CO/0204005 (2002).

[RSZ] A. Robertson, D. Saracino, and D. Zeilberger,
Refined restricted permutations,\\
arXiv:math.CO/0203022 (2002).

[S] C. Schensted, Longest increasing and decreasing subsequences,
{\it Canad. J. Math.} {\bf 13} (1961), 179-191.

[SiS] R. Simion and F. Schmidt, Restricted permutations,
{\it European Journal of Combinatorics} {\bf 6} (1985), 383-406.

[Z] D. Zeilberger, Rodica Simion (1955-2000):  An  (almost) perfect
enumerator and human being, {\it Shalosh B. Ekhad's and 
Doron Zeilberger's Personal Journal}, {\tt
http://www.math.rutgers.edu/$\sim$zeilberg/pj.html}.

\newpage
\section*{\normalsize Appendix}
Below we provide values of $i_n^k(\alpha)$ and
$i_n^k(\emptyset;\alpha)$ for small $n$
and all $\alpha \in S_3$. 
 
\vskip 20pt
\footnotesize 
\setlength{\arraycolsep}{2.5pt}
\renewcommand{\arraystretch}{.85}
$$
\begin{array}{lllll} 
\begin{array}{lcccccccccccc} 
{}_n \diagdown^k \!\! &\vline&^0&^1&^2&^3&^4&^5&^6&^7&^8 \\
\hline 
_0&\vline&1\\
_1&\vline&0&1\\
_2&\vline&1&0&1\\
_3&\vline&0&3&0&0\\
_4&\vline&3&0&3&0&0\\
_5&\vline&0&10&0&0&0&0\\
_6&\vline&10&0&10&0&0&0&0\\
_7&\vline&0&35&0&0&0&0&0&0\\
_8&\vline&35&0&35&0&0&0&0&0&0\\ 
\end{array}
&\hskip 10pt
&
\begin{array}{lcccccccccccc} 
{}_n \diagdown^k \!\! &\vline&^0&^1&^2&^3&^4&^5&^6&^7&^8 \\
\hline 
_0&\vline&1\\
_1&\vline&0&1\\
_2&\vline&1&0&1\\
_3&\vline&0&2&0&1\\
_4&\vline&2&0&3&0&1\\
_5&\vline&0&5&0&4&0&1\\
_6&\vline&5&0&9&0&5&0&1\\
_7&\vline&0&14&0&14&0&6&0&1\\
_8&\vline&14&0&28&0&20&0&7&0&1\\
\end{array}
&\hskip 10pt
&
\begin{array}{lcccccccccccc} 
{}_n \diagdown^k \!\! &\vline&^0&^1&^2&^3&^4&^5&^6&^7&^8 \\
\hline 
_0&\vline&1\\
_1&\vline&0&1\\
_2&\vline&1&0&1\\
_3&\vline&0&3&0&1\\
_4&\vline&2&0&5&0&1\\
_5&\vline&0&8&0&7&0&1\\
_6&\vline&4&0&18&0&9&0&1\\
_7&\vline&0&20&0&32&0&11&0&1\\
_8&\vline&8&0&56&0&50&0&13&0&1\\
\end{array}
\\ \\
\mathbf{i_n^k(123)}&&\mathbf{i_n^k(132)=i_n^k(321)=i_n^k(213)}
&&\mathbf{i_n^k(231)=i_n^k(312)}\\
\\
\\
\begin{array}{lcccccccccccc} 
{}_n \diagdown^k \!\! &\vline&^0&^1&^2&^3&^4&^5&^6&^7&^8 \\
\hline 
_0&\vline&0\\
_1&\vline&0&0\\
_2&\vline&0&0&0\\
_3&\vline&0&0&0&1\\
_4&\vline&0&0&0&0&0\\
_5&\vline&0&0&0&3&0&0\\
_6&\vline&0&0&0&0&0&0&0\\
_7&\vline&0&0&0&9&0&0&0&0\\
_8&\vline&0&0&0&0&0&0&0&0&0\\
\end{array}
&\hskip 10pt
&
\begin{array}{lcccccccccccc} 
{}_n \diagdown^k \!\! &\vline&^0&^1&^2&^3&^4&^5&^6&^7&^8 \\
\hline 
_0&\vline&0\\
_1&\vline&0&0\\
_2&\vline&0&0&0\\
_3&\vline&0&1&0&0\\
_4&\vline&0&0&1&0&0\\
_5&\vline&0&2&0&1&0&0\\
_6&\vline&0&0&3&0&1&0&0\\
_7&\vline&0&5&0&4&0&1&0&0\\
_8&\vline&0&0&9&0&5&0&1&0&0\\
\end{array}
&\hskip 10pt
&
\begin{array}{lcccccccccccc} 
{}_n \diagdown^k \!\! &\vline&^0&^1&^2&^3&^4&^5&^6&^7&^8 \\
\hline 
_0&\vline&0\\
_1&\vline&0&0\\
_2&\vline&0&0&0\\
_3&\vline&0&0&0&0\\
_4&\vline&0&0&1&0&0\\
_5&\vline&0&0&0&2&0&0\\
_6&\vline&0&0&2&0&3&0&0\\
_7&\vline&0&0&0&8&0&4&0&0\\
_8&\vline&0&0&5&0&18&0&5&0&0\\
\end{array}
\\
\\
\mathbf{i_n^k(\emptyset;123)}&&
\mathbf{i_n^k(\emptyset;132)=i_n^k(\emptyset;213)}&&
\mathbf{i_n^k(\emptyset;231)=i_n^k(\emptyset;312)}\\
\\
\\
\begin{array}{lcccccccccccc} 
{}_n \diagdown^k \!\! &\vline&^0&^1&^2&^3&^4&^5&^6&^7&^8 \\
\hline 
_0&\vline&0\\
_1&\vline&0&0\\
_2&\vline&0&0&0\\
_3&\vline&0&1&0&0\\
_4&\vline&0&0&2&0&0\\
_5&\vline&0&4&0&3&0&0\\
_6&\vline&0&0&10&0&4&0&0\\
_7&\vline&0&14&0&18&0&5&0&0\\
_8&\vline&0&0&40&0&28&0&6&0&0\\
\end{array}
\\ \\
\mathbf{i_n^k(\emptyset;321)}
\end{array}
$$

\setlength{\arraycolsep}{5pt}
\renewcommand{\arraystretch}{1.0}

\end{document}